\documentclass{amsart}
\theoremstyle{plain}
\newtheorem{Thm}{Theorem}

\errorcontextlines=0 \numberwithin{equation}{section}

\begin{document}

\title[gradient estimate]
{Gradient estimate for the Poisson equation and the non-homogeneous
heat equation on compact Riemannian manifolds}

\author{Li Ma, Liang Cheng}

\address{Department of mathematical sciences \\
Tsinghua university \\
Beijing 100084 \\
China} \email{lma@math.tsinghua.edu.cn} \dedicatory{}
\date{May 26th, 2009}

\begin{abstract}

In this short note, we study the  gradient estimate of positive
solutions to Poisson equation and the non-homogeneous heat
equation in a compact Riemannian manifold $(M^n,g)$. Our results
extend the gradient estimate for positive harmonic functions and
positive solutions to heat equations.


\end{abstract}

\thanks{$^*$ The research is partially supported by the National Natural Science
Foundation of China 10631020 and SRFDP 20060003002. }
 \maketitle

\section{Introduction}
In the process of the study of the positive solutions to non-local
non-homogeneous heat equation (see \cite{CL09},\cite{MA}, \cite{MC},
and \cite{MC2}) in a compact Riemannian manifold $(M,g)$:
$$
u_t=\Delta u+\lambda(t)u+f(x,t),
$$
with the initial data  $u(0,x)=u_0(x)$, where $\int_Mu_0(x)^2=1$ and
$$\lambda (t)=\int_M|\nabla u|^2-fu,
$$
we find that it is interesting to study the gradient estimate for
positive solutions to the elliptic equation
$$
-\Delta u=A(x), \; \; in \; M
$$
and the heat equation
$$
(\partial_t-\Delta) u=A(x,t), \; \;in \; M\times (0,T).
$$
We shall follow the ideas in \cite{LM} and \cite{ML}, which uses
tricks from the works of Cheng-Yau \cite{SY} on harmonic functions
and Li-Yau \cite{SY} on heat equations. However, some new
contributions have to be provided since they have treated different
situation as ours (see \cite{BC},\cite{Evans},\cite{LM} and
\cite{P02}). It is also clear that our result can be extended to
complete Riemannian manifolds. For more related works on complete
manifolds, one may look at \cite{Evans}, \cite{H95}, \cite{DM}, and
\cite{LCZ}.

It is a surprise to us that there is few literature about the
gradient estimate for positive solutions to Poisson equations and
non-homogeneous heat equations. With this understanding and
motivated from the work of L.Caffarelli and F.H.Lin \cite{CL09}
and our paper \cite{LM}, we consider the gradient estimates for
Poisson equation and non-homogeneous heat equation in a compact
Riemannian manifold.

 We shall always assume that $Ric(g)\geq K$ on $M$ for some real
constant $K$.

For the Poisson equation, we have the following result.

\begin{Thm}\label{thm1}Let $u>0$ be a smooth solution to the Poisson equation on $M$
$$
-\Delta u=A(x).
$$
Then we have
$$
|\nabla w|^2+A(x)u^{-1}\leq 2n \sup\{K-Au^{-1},
A^2u^{-2}-[n^{-1}4(Au^{-1}-K)^2+2KAu^{-1}+u^{-1}\Delta A]\},
$$
\end{Thm}

For the non-homogeneous heat equation on $M\times [0,T)$, we have
the following result.

\begin{Thm}\label{thm2}
Let $u>0$ be a smooth solution to the non-homogeneous heat
equation on $M\times [0,T)$
$$
(\partial_t-\Delta) u=A(x,t).
$$
Let, for $a>1$, $$ F=t(|\nabla w|^2+aA(x,t)u^{-1}-aw_t). $$ Then
there is a constant $C(u^{-1},|A|,|\nabla A|,|\Delta A|,K,a,T)>0$
such that
$$\sup_{M\times(0,T)}F\leq C(u^{-1},|A|,|\nabla A|,|\Delta A|,K,a,T).$$
\end{Thm}
Related local gradient estimates can be extended to complete
non-compact Riemannian manifolds, which will appear elsewhere.

 This paper is organized as
follows. In section \ref{sect2} we prove Theorem \ref{thm1}. In
section \ref{sect3} we do the gradient estimate for a positive
smooth solution to the non-homogeneous heat equation.

\section{gradient estimate for Poisson equation}\label{sect2}
We firstly prove Theorem \ref{thm1} about the gradient estimate for
parabolic equations. Recall here that we are considering the object
for Poisson equation on  the compact Riemannian manifold $(M^n,g)$.

We now recall the famous Bochner formula any smooth function $v$ on
a Riemannian manifold $(M^n,g)$:
\begin{equation}\label{bochner}
\Delta |\nabla v|^2=2|D^2v|^2+2(\nabla v,\nabla\Delta v)+2Ric(\nabla
v,\nabla v).
\end{equation}
Recall that $|D^2v|^2\geq \frac{1}{n}|\Delta v|^2$. So we have
$$
\Delta |\nabla v|^2\geq\frac{2}{n}|\Delta v|^2+2(\nabla
v,\nabla\Delta v)+2Ric(\nabla v,\nabla v)
$$
This formula will plays a key role in our gradient estimate.

Let $u>0$ be a smooth solution to the Poisson equation on $M$
$$
-\Delta u=A(x).
$$
Set
$$
w=log u.
$$
Then we have
\begin{equation}\label{ediff}
-\Delta w=|\nabla w|^2+A(x)u^{-1}.
\end{equation}

Let $Q=|\nabla w|^2+A(x)u^{-1}$ be the Harnack quantity. Then
$Q=-\Delta w$.

 By (\ref{ediff}), we obtain that
$$
\Delta Q=\Delta |\nabla w|^2+\Delta (A(x)u^{-1}).
$$
 Using the Bochner formula
(\ref{bochner}), we get
$$
\Delta |\nabla w|^2\geq \frac{2}{n}Q^2+2(\nabla v,\nabla
Q)-2K|\nabla w|^2.
$$
Then we have
$$
\Delta Q\geq \frac{2}{n}Q^2+2(\nabla w,\nabla
Q)-2K(Q-A(x)u^{-1})+\Delta (A(x)u^{-1})
$$
 Note that
$$
\Delta (A(x)u^{-1})=u^{-1}\Delta A(x)-2u^{-1}\nabla A\cdot\nabla
w+Au^{-1}(2Q-Au^{-1}).
$$
Then at the maximum point $p\in M$ of $Q$ (which can be assumed
positive), we have
$$ \Delta Q\leq 0, \; \nabla Q=0.
$$
Then we have
$$
0\geq \frac{2}{n}Q^2+(2Au^{-1}-2K)Q+2KAu^{-1}+u^{-1}\Delta
A-A^2u^{-2}-2u^{-1}\nabla A\cdot\nabla w.
$$
Using the Cauchy-Schwartz inequality we get for any $b>0$,
$$
\frac{2}{n}[Q+n(Au^{-1}-b-K)]^2+n^{-1}4(Au^{-1}-b-K)^2+2KAu^{-1}+u^{-2}(u\Delta
A-\frac{|\nabla A|^2}{2b}-A^2)\leq 0.
$$

If $Q> 2n(K+b-Au^{-1})$, then
$$
Q+n(Au^{-1}-b-K)> Q/2>0.
$$
Hence we have
$$
\frac{1}{2n}Q^2\leq (A^2+\frac{|\nabla
A|^2}{2b})u^{-2}-[n^{-1}4(Au^{-1}-K)^2+2KAu^{-1}+u^{-1}\Delta A].
$$
In conclusion we have
$$
Q\leq 2n \sup\{K+b-Au^{-1},(A^2+\frac{|\nabla A|^2}{2b})u^{-2}
-[n^{-1}4(Au^{-1}-K)^2+2KAu^{-1}+u^{-1}\Delta A]\}
$$

This implies that by choosing $b=1/2$,
$$
|\nabla w|^2+A(x)u^{-1}\leq 2n \sup\{K+\frac{1}{2}-Au^{-1},
(A^2+|\nabla A|^2)u^{-2}-[n^{-1}4(Au^{-1}-K)^2+u^{-1}(2KA+\Delta
A)]\},
$$
which is the gradient estimate wanted for positive solutions to the
Poisson equation. This completes the proof of Theorem \ref{thm1}.

\section{gradient estimate for non-homogeneous heat equation}\label{sect3}
We now prove Theorem \ref{thm2}. Let $u>0$ be a smooth solution to
the non-homogeneous heat equation on $M\times [0,T)$
\begin{equation}\label{heat}
(\partial_t-\Delta) u=A(x,t).
\end{equation}
Set
$$
w=log u.
$$
Then we have
\begin{equation}\label{ediff}
(\partial_t-\Delta) w=|\nabla w|^2+Au^{-1}.
\end{equation}

Following Li-Yau \cite{SY} we let $F=t(|\nabla w|^2+aAu^{-1}-aw_t)$
(where $a>1$) be the Harnack quantity for (\ref{heat}). Then we have
$$ |\nabla w|^2=\frac{F}{t}-aAu^{-1}+aw_t,
$$
$$
\Delta w=w_t-|\nabla
w|^2-Au^{-1}=-\frac{F}{at}-(1-\frac{1}{a})|\nabla w|^2.
$$
and
$$
w_t-\Delta w=|\nabla w|^2+Au^{-1}=\frac{F}{t}+(1-a)Au^{-1}+aw_t.
$$
 Note that
$$
(\partial_t-\Delta) w_t=2\nabla w\nabla w_t+\frac{d}{dt}(Au^{-1}).
$$
Using the Bochner formula, we have
 $$ (\partial_t-\Delta) |\nabla
w|^2=2\nabla w\nabla w_t-[2|D^2w|^2+2(\nabla w,\nabla\Delta
w)+2Ric(\nabla w,\nabla w)],
$$
and using (\ref{ediff}) we get
$$
(\partial_t-\Delta) |\nabla w|^2=2\nabla w\nabla (w_t-\Delta
w)-[2|D^2w|^2+2Ric(\nabla w,\nabla w)],
$$
which can be rewritten as
$$
(\partial_t-\Delta) |\nabla w|^2=2\nabla w\nabla
[\frac{F}{t}+(1-a)Au^{-1}+aw_t]-[2|D^2w|^2+2Ric(\nabla w,\nabla w)].
$$
Then we have $$ (\partial_t-\Delta) (|\nabla w|^2-aw_t) =2\nabla
w\nabla [\frac{F}{t}+(1-a)Au^{-1}]
$$
$$-[2|D^2w|^2+2Ric(\nabla w,\nabla
w)]-a\frac{d}{dt}(Au^{-1}).
$$
Hence
\begin{eqnarray*}
&&(\partial_t-\Delta) (|\nabla w|^2-aw_t+aAu^{-1})\\
&=&(\partial_t-\Delta) (|\nabla w|^2-aw_t)+a(\partial_t-\Delta)
Au^{-1})\\
&=& 2\nabla w\nabla [\frac{F}{t}+(1-a)Au^{-1}]
-[2|D^2w|^2+2Ric(\nabla w,\nabla
w)]-a\frac{d}{dt}(Au^{-1})\\
& &+a(\partial_t-\Delta) Au^{-1})\\
&=&2\nabla w\nabla [\frac{F}{t}+(1-a)(Au^{-1})]
-[2|D^2w|^2+2Ric(\nabla w,\nabla w)]-a\Delta (Au^{-1}).
\end{eqnarray*}
Then we have
$$
(\partial_t-\Delta)F =\frac{F}{t}+2t\nabla w\nabla
[\frac{F}{t}+(1-a)(Au^{-1})]
$$
$$-t[2|D^2w|^2+2Ric(\nabla w,\nabla
w)]-at\Delta (Au^{-1}).
$$
Assume that $$ \sup_{M\times[0,T]} F>0.
$$
Applying the maximum principle at the maximum point $(z,s)$, we then
have
$$ (\partial_t-\Delta)F\geq 0, \; \nabla F=0.
$$
In the following our computation is always at the point $(z,s)$. So
we get
\begin{equation}\label{key1}
\frac{F}{s}+2(1-a)s\nabla w\nabla (Au^{-1})-s[2|D^2w|^2+2Ric(\nabla
w,\nabla w)]-as\Delta (Au^{-1})\geq 0.
\end{equation}
That is
\begin{equation}\label{key2}
F-as^2\Delta (Au^{-1})\geq 2(a-1)s^2\nabla w\nabla (Au^{-1})+
s^2[2|D^2w|^2+2Ric(\nabla w,\nabla w)].
\end{equation}
Set
$$
\mu=\frac{|\nabla w|^2}{F}{(z,s)}.
$$
Then at $(z,s)$,$$ |\nabla w|^2=\mu F.
$$
Hence
$$
\nabla \frac{A}{u}=\frac{\nabla A}{u}-\frac{A\nabla
u}{u^2}=\frac{\nabla A}{u}-\frac{A}{u}\nabla w.
$$
So
$$
\nabla w \cdot \nabla \frac{A}{u}=\frac{\nabla w \cdot\nabla
A}{u}-\frac{A}{u}|\nabla w|^2\geq -\frac{|\nabla w||\nabla
A|}{u}-\frac{A}{u}|\nabla w|^2
$$
\begin{equation}\label{simplify1}
=-\frac{|\nabla A|}{u}\sqrt{\mu F}-\frac{A}{u}\mu F \geq
-\frac{1}{2}\frac{|\nabla A|^2}{u}-(\frac{1}{2}+A)\frac{\mu F}{u}.
\end{equation}
Further more, we have
\begin{eqnarray*}
(\partial_t-\Delta)(Au^{-1})&=&\frac{1}{u}(\partial_t-\Delta)A-\frac{A}{u^2}(\partial_t-\Delta)u+
\frac{2}{u^2}\nabla u \cdot \nabla A-2\frac{A}{u^3}|\nabla u|^2\\
&=&\frac{1}{u}(\partial_t-\Delta)A-\frac{A^2}{u^2}+
\frac{2}{u}\nabla w \cdot \nabla A-2\frac{A}{u}|\nabla w|^2\\
&\leq&\frac{1}{u}(\partial_t-\Delta)A-\frac{A^2}{u^2}+
\frac{2}{u}\sqrt{\mu F}|\nabla A|-2\frac{A}{u}\mu F\\
&\leq&\frac{1}{u}(\partial_t-\Delta)A-\frac{A^2}{u^2}+ \frac{\mu
F}{u}+\frac{|\nabla A|^2}{u}-2\frac{A}{u}\mu F,
\end{eqnarray*}
and
\begin{eqnarray*}
\partial_t(Au^{-1})&=&\frac{A_t}{u}-\frac{A}{u^2}u_t\\
&=&\frac{A_t}{u}-\frac{A}{u}w_t\\
&=&\frac{A_t}{u}-\frac{A}{u}(\frac{1}{a}(|\nabla
w|^2-\frac{F}{s})+Au^{-1})\\
&=&\frac{A_t}{u}-\frac{A}{u}\cdot
\frac{F}{a}(\mu-\frac{1}{s})-\frac{A^2}{u^2}.
\end{eqnarray*}
Hence
\begin{equation}\label{simplify2}
-\Delta(Au^{-1})=(\partial_t-\Delta)(Au^{-1})-\partial_t(Au^{-1})
\end{equation}
$$
\leq-\frac{\Delta A}{u}+\frac{|\nabla A|^2}{u}+\frac{1-2A}{u}\mu
F+\frac{A}{u}\cdot \frac{F}{a}(\mu-\frac{1}{s})
$$
$$
<-\frac{\Delta A}{u}+\frac{|\nabla A|^2}{u}+\frac{1-2A}{u}\mu
F+\frac{A}{u}\cdot \frac{F}{a}\mu. \ \ \ \ \ \ \ \ \
$$
Note that
$$
|D^2w|^2+Ric(\nabla w,\nabla w)\geq \frac{1}{n}|\Delta w|^2-K|\nabla
w|^2.
$$
So
$$
|D^2w|^2+Ric(\nabla w,\nabla w)\geq
\frac{1}{n}(\frac{F}{as}+(1-\frac{1}{a})|\nabla w|^2)^2-K|\nabla
w|^2
$$
\begin{equation}\label{simplify3}
=\frac{F^2}{n}(\frac{1}{as}+(1-\frac{1}{a})\mu)^2-K\mu F.
\end{equation}
Substitute (\ref{simplify1}) (\ref{simplify2}) and (\ref{simplify3})
into (\ref{key2}), we get
$$
F+\frac{as^2}{u}(-\Delta A+|\nabla A|^2)+\mu F
\frac{s^2}{u}(a+(1-2a)A)
$$
$$
\geq -s^2(a-1)\frac{|\nabla A|^2}{u}-(a-1)s^2\frac{1+2A}{u}\mu F
+\frac{2F^2}{n}(\frac{1}{a}+(1-\frac{1}{a})\mu s)^2-2s^2K\mu F.
$$
Assume that
$$
F\geq \frac{as^2}{u}(-\Delta A+|\nabla A|^2)+ s^2(a-1)\frac{|\nabla
A|^2}{u},
$$
for otherwise we are done. Then we have
$$
2F+\mu F \frac{s^2}{u}(a+(1-2a)A)+(a-1)s^2\frac{1+2A}{u}\mu
F+2s^2K\mu F
$$
$$
\geq\frac{2F^2}{n}(\frac{1}{a}+(1-\frac{1}{a})\mu s)^2.
$$

Simplify this inequality, we get
$$
\frac{2F}{n}\frac{1}{a^2}\leq \frac{2}{(1+(a-1)\mu s)^2}+\frac{\mu
s}{(1+(a-1)\mu s)^2}
$$
$$
\cdot s(u^{-1}(a+(1-2a)A)+u^{-1}(a-1)(1+2A)+2K).
$$
Hence we have the estimate for $F$ at $(z,s)$ such that
$$
F(z,s)\leq C(u^{-1},|A|,|\nabla A|,|\Delta A|,K,a,T),
$$
which is the desired gradient estimate. This completes the proof
of Theorem \ref{thm2}.


\begin{thebibliography}{20}
\bibitem{A98}
T. Aubin, \emph{Some Nonlinear Problems in Riemannian Geometry},
Springer Monogr. Math., Springer-Verlag, Berlin, 1998.

\bibitem{CL09}
C.Caffarelli, F.Lin, \emph{Nonlocal heat flows preserving the
$L^2$ energy}, Discrete and continuous dynamical systems. 23,
49-64 (2009).

\bibitem{cdy} K.C.Chang, W.Y.Ding, and R.Ye, \emph{Finite time blow-up
of the heat flow of harmonic maps from surfaces}, JDG,
36(1992)507-515.

\bibitem{BC} B.Chow, P.Lu,L.Ni. \emph{Hamilton's Ricci Flow}. Science
Press. American Mathematical Society, Beijing.Providence(2006).

\bibitem{BH97}
Ben Chow and Richard Hamilton, \emph{Constrained and linear Harnack
inequalities for parabolic equations}, Inventiones Mathematicae 129,
213-238 (1997).

\bibitem{DM} Xianzhe Dai and Li Ma, \emph{Mass under Ricci flow},
Commun. Math. Phys., 274, 65-80 (2007).

\bibitem{Evans} L.Evans,
\emph{Partial Differential Equations,} Graduate studies in Math.,
AMS, 1986

\bibitem{H95}
 R.Hamilton, {\em The formation of Singularities in the Ricci flow},
 Surveys in Diff. Geom.,
 Vol.2, pp7-136, 1995.

 \bibitem{LM}
L.Ma, \emph{Gradient estimates for a simple elliptic equation on
complete non-compact Riemannian manifolds}, Journal of Functional
Analysis, 241(2006)374-382.

 \bibitem{MA} L.Ma, A.Q.Zhu, \emph {On a length preserving curve
flow.} Preprint,2008.

\bibitem{MC} L.Ma, L.Cheng, \emph {A non-local area preserving curve
flow}, Preprint, 2008.

\bibitem{MC2} L.Ma, L.Cheng, \emph {non-local heat flows
and gradient estimates on closed manifolds}, Preprint, 2009.

\bibitem{ML}
Li Ma, Baiyu, Liu, \emph{Convex eigenfunction of a drifting
Laplacian operator and the fundamental gap}, Pacific Journal of
Math., 240(2009)343-361

\bibitem{LCZ} L.Ma, Chong Li, and Lin Zhao, \emph{Monotone solutions to a class of
elliptic and diffusion equations}, CPAA, 6(2007)237-246.

\bibitem{P02} Grisha Perelman,
\emph{The entropy formula for the Ricci flow and its geometric
applications},math.DG/0211159,2002.

\bibitem{SY2} R.Schoen,
\emph{Analytic aspects for Harmonic maps}, Seminar in PDE, edited by
S.S.Chern, Springer, 1984.

\bibitem{SY} R.Schoen and S.T.Yau,
\emph{Lectures on Differential Geometry}, international Press, 1994.
\end{thebibliography}
\end{document}